\documentclass[preprint,12pt]{elsarticle}

\usepackage{amssymb}
\usepackage{mathrsfs}
\usepackage{amsmath}
\usepackage{float}
\usepackage{amsthm,amsfonts}
\usepackage{graphicx}
\newtheorem{lemma}{Lemma}[section]
\newtheorem{theorem}[lemma]{Theorem}
\newtheorem{corollary}[lemma]{Corollary}
\newtheorem{problem}[lemma]{Problem}

\newenvironment{prof}[1][Proof]{\noindent\textit{#1}\quad }

\journal{}

\begin{document}

\newcommand{\js}{\hfill $\square$}

\begin{frontmatter}



\title{Size of edge-critical uniquely 3-colorable planar graphs
\footnote{Supported by 973 Program of China 2013CB329601, 2013CB329603, National Natural
Science Foundation of China Grant 61309015 and National Natural
Science Foundation of China Special Equipment Grant 61127005.
}}


%\author{}

%\address{}

\author[PKU]{Zepeng Li}\ead{lizepeng@pku.edu.cn}
\author[PKU]{Enqiang Zhu}\ead{zhuenqiang@pku.edu.cn}
\author[CDU1,CDU2]{Zehui Shao}\ead{zshao@cdu.edu.cn}
\author[PKU]{Jin Xu}\ead{jxu@pku.edu.cn}
\address[PKU]{ Key Laboratory of High Confidence Software Technologies, Peking University, Beijing, 100871, China}
\address[CDU1]{Key Laboratory of Pattern Recognition and Intelligent Information Processing, Institutions of Higher Education of Sichuan Province, China}
\address[CDU2]{School of Information Science and Technology, Chengdu University,  Chengdu, 610106,  China}

\begin{abstract}
A graph $G$ is \emph{uniquely k-colorable} if the chromatic number of $G$ is $k$ and $G$ has
only one $k$-coloring up to permutation of the colors. A uniquely $k$-colorable graph $G$
is edge-critical if $G-e$ is not a uniquely $k$-colorable graph for any edge $e\in E(G)$. Mel'nikov and Steinberg [L. S. Mel'nikov, R. Steinberg, One counterexample for two conjectures on three coloring, Discrete Math. 20 (1977) 203-206] asked to find an exact upper bound for the number of edges in a edge-critical 3-colorable planar graph with $n$ vertices.
In this paper, we give some properties of edge-critical uniquely 3-colorable planar graphs and prove that if $G$ is such a graph with $n(\geq6)$ vertices, then $|E(G)|\leq \frac{5}{2}n-6 $, which improves the upper bound $\frac{8}{3}n-\frac{17}{3}$ given by Matsumoto [N. Matsumoto, The size of edge-critical uniquely 3-colorable planar graphs, Electron. J. Combin. 20 (3) (2013) $\#$P49]. Furthermore, we find some edge-critical 3-colorable planar graphs which have $n(=10,12, 14)$ vertices and $\frac{5}{2}n-7$ edges.
\end{abstract}

\begin{keyword}
planar graph; unique coloring; uniquely $3$-colorable planar graph; edge-critical




\MSC 05C15
\end{keyword}

\end{frontmatter}


\section{Introduction}

%plantri -p c2m3 9 9.g6





 A graph $G$ is \emph{uniquely k-colorable} if $\chi(G)=k$ and $G$ has only one $k$-coloring
up to permutation of the colors, where the coloring is called a unique $k$-coloring. In other words, all $k$-colorings
of $G$ induce the same partition of $V(G)$ into $k$ independent sets. In addition, uniquely colorable graphs may be defined in terms of their chromatic polynomials, which initiated by Birkhoff \cite{Birkhoff1912} for planar graphs in 1912 and, for general graphs, by Whitney \cite{Whitney1932} in 1932. Because a graph $G$ is uniquely $k$-colorable if and only if its chromatic polynomial is $k!$. For a discussion of chromatic polynomials, see Read \cite{Read1968}.

Let $G$ be a uniquely $k$-colorable graph, $G$ is \emph{edge-critical} if $G-e$ is not uniquely $k$-colorable for any edge $e\in E(G)$.
Uniquely colorable graphs were defined and studied firstly by Harary and Cartwright \cite{Harary1968} in 1968. They proved the following theorem.

\begin{theorem}(Harary and Cartwright \cite{Harary1968}) \label{theorem1.1}
Let $G$ be a uniquely $k$-colorable graph. Then for any unique $k$-coloring of $G$, the subgraph induced by the union of any two color classes is connected.
\end{theorem}

As a corollary of Theorem \ref{theorem1.1}, it can be seen that a uniquely $k$-colorable graph $G$ has at least $(k-1)|V(G)|-{k \choose 2}$ edges. Furthermore, if a  uniquely $k$-colorable graph $G$ has exactly $(k-1)|V(G)|-{k \choose 2}$ edges, then $G$ is edge-critical. There are many references on uniquely colorable graphs. For example see Chartrand and Geller \cite{Chartrand1969},  Harary, Hedetniemi and Robinson \cite{Harary1969} and
Bollob\'{a}s \cite{Bollob¨¢s1978}.

Chartrand and Geller \cite{Chartrand1969} in 1969 started to study uniquely colorable planar graphs. They proved that uniquely 3-colorable planar graphs with at least 4 vertices contain at least two triangles, uniquely 4-colorable planar graphs are maximal planar graphs, and uniquely 5-colorable planar graphs do not exist. Aksionov \cite{Aksionov1977} in 1977 improved the low bound for the number of triangles in a uniquely 3-colorable planar graph. He proved that a uniquely 3-colorable planar graph with at least 5 vertices contains at least 3 triangles and gave a complete description of uniquely 3-colorable planar graphs containing exactly 3 triangles.

For an edge-critical uniquely $k$-colorable planar graph $G$, if $k=2$, then it is easy to deduce that $G$ is tree and has exactly $|V(G)|-1$ edges. If $k=4$, then $G$ is a maximal planar graph and has exactly $3|V(G)|-6$ edges by Euler's Formula. Therefore, it is sufficient to consider the size of uniquely $3$-colorable planar graphs. We denote by $\mathcal{U}_E$ the set of all edge-critical uniquely $3$-colorable planar graphs and by $size(n)$ the upper bound of the size of edge-critical uniquely $3$-colorable planar graphs with $n$ vertices.

In 1977 Aksionov \cite{Aksionov1977} conjectured that $size(n)=2n-3$. However, in the same year,
Mel'nikov and Steinberg \cite{Mel'nikov1977} disproved the conjecture by constructing a counterexample $H$, which has 16 vertices and 30 edges. Moreover, they proposed the following problems:

\begin{problem}(Mel'nikov and Steinberg \cite{Mel'nikov1977})\label{problem1}
Find an exact upper bound for the number of edges in a edge-critical 3-colorable planar graph with $n$ vertices. Is it true that $size(n)=\frac{9}{4}n-6$ for any $n\geq 12$?
\end{problem}

Recently, Matsumoto \cite{Matsumoto2013} constructed an infinite family of
edge-critical uniquely 3-colorable planar graphs with $n$ vertices and $\frac{9}{4}n-6$ edges, where $n\equiv 0 (\textrm{mod}~4)$. He also gave a non-trivial upper bound $\frac{8}{3}n-\frac{17}{3}$ for $size(n)$.

In this paper, we give some properties of edge-critical uniquely 3-colorable planar graphs with $n$ vertices and improve the upper bound of $size(n)$ given by Matsumoto \cite{Matsumoto2013} to $\frac{5}{2}n-6$, where $n\geq 6$. Moreover, we give some edge-critical 3-colorable planar graphs which have $n(=10,12, 14)$ vertices and $\frac{5}{2}n-7$ edges. It follows that the conjecture of Mel'nikov and Steinberg \cite{Mel'nikov1977} is false because $\frac{5}{2}n-7> \frac{9}{4}n-6$ if $n\geq 12$.

\section{Notation}

Only finite, undirected and simple graphs are considered in this paper. For a planar graph $G=(V(G), E(G),F(G))$, $V(G)$, $E(G)$ and $F(G)$ are the sets of vertices, edges and faces of $G$, respectively. We denote by $\delta(G)$ and $\Delta(G)$  the \emph{minimum degree} and \emph{maximum degree} of graph $G$. The degree of a vertex $v \in V(G)$, denoted by $d_{G}(v)$, is the number of neighbors of $v$ in $G$. The degree of a face $f\in F(G)$, denoted by $d_G(f)$, is the number of edges in its boundary, cut edges being counted twice. When no confusion can arise, $d_{G}(v)$ and $d_{G}(f)$ are simplified by $d(v)$ and $d(f)$, respectively. A face $f$ is a \emph{$k$-face} if $d_G(f)=k$ and a \emph{$\geq$k-face} if $d_G(f)\geq k$. The similar notation is used for cycles.  We denote by $V_i(G)$ the set of vertices of $G$ with degree $i$ and by $V_{\geq i}(G)$ the set of vertices of $G$ with degree at least $i$, where $\delta(G)\leq i \leq \Delta(G)$.  The similar notation is used for the set of faces of $G$.

A $k$-wheel is the graph consists of a single vertex $v$ and a cycle $C$ with $k$ vertices together with $k$ edges from $v$ to each vertex of $C$. A planar (resp. outerplanar) graph $G$ is \emph{maximal} if $G+uv$ is not planar (resp. outerplanar) for any two nonadjacent vertices $u$ and $v$ of $G$. Let $V_1$ and $V_2$ be two disjoint subset of $V(G)$, we use $e(V_1,V_2)$ to denote the number of edges of $G$ with one end in $V_1$ and the other in $V_2$. In particular, if $V_1$ or $V_2=\{v\}$, we simply write $e(v,V_2)$ or $e(V_1,v)$ for $e(V_1,V_2)$, respectively.
To \emph{contract} an edge $e$ of a graph $G$ is to delete the edge and then identify its ends. The resulting graph is denoted by $G/e$.
Two faces $f_1$ and $f_2$ of $G$ are \emph{adjacent} if they have at least one common edge. A $k$-cycle $C$ is said to be a \emph{separating} $k$-\emph{cycle} in $G$ if the removal of $C$ disconnects the graph $G$.

A \emph{k-coloring} of $G$ is an assignment of $k$ colors to $V(G)$ such that no two adjacent vertices are assigned the same color. Naturally, a $k$-coloring can be viewed as a partition $\{V_1,V_2,\cdots,V_k\}$ of $V$, where $V_i$ denotes the set of vertices assigned color $i$, and is called a \emph{color class} of the coloring for any $i=1,2,\cdots,k$. Two $k$-colorings $f$ and $f^\prime$ of $G$ are said to be \emph{distinct}
if they produce two distinct partitions of $V(G)$ into $k$ color classes. A graph $G$ is \emph{k-colorable} if there exists a $k$-coloring of $G$, and
the \emph{chromatic number} of $G$, denoted by $\chi(G)$, is the minimum number $k$ such that $G$ is
$k$-colorable.

The notations and terminologies not mentioned here can be found in \cite{Bondy2008}.

\section{Properties of edge-critical uniquely $3$-colorable planar graphs}

Let $G$ be a 3-colorable planar graph and $f$ be a 3-coloring of $G$. It is easy to see that the restriction of $f$ to $G-e$ is a 3-coloring of $G-e$, where $e\in E(G)$. For convenience, we also say $f$ is a 3-coloring of $G-e$. If there exists a 3-coloring $f'$ of $G-uv$ such that $f'(u)\neq f'(v)$, then we say that $f'$ can be \emph{extended} to a 3-coloring of $G$.

\begin{theorem}\label{theorem2.1}
Let $G$ be a uniquely 3-colorable planar graph. Then $G \in \mathcal{U}_E$ if and only if $G/ e$ is 3-colorable for any edge $e\in E(G)$.
\end{theorem}
\begin{prof}
Suppose that $G \in \mathcal{U}_E$, then, by definition, $G-e$ has at least two distinct 3-colorings for each $e=uv\in E(G)$. Since $G$ is uniquely 3-colorable, we conclude that there exists a 3-coloring $f$ of $G-e$ such that $f(u)=f(v)$. Hence $G/ e$ is 3-colorable.

Conversely, suppose that $G \notin \mathcal{U}_E$. Then there exists an edge $e'=uv\in E(G)$ such that $G-e'$ is also a uniquely 3-colorable planar graph. Obviously, for any unique 3-coloring $f$ of $G$, we have $f(u)\neq f(v)$. So $G/ e'$ is not 3-colorable. This establishes Theorem \ref{theorem2.1}. \qed
\end{prof}

The following result is obtained by Theorem \ref{theorem2.1}.

\begin{corollary}\label{corollary2.2}
Let $G \in \mathcal{U}_E$ and $v\in V(G)$. If $v$ is incident with exactly one 4-face and all other faces
incident with $v$ are triangular, then $d(v)$ is even.
\end{corollary}
\begin{prof}
Suppose that the result is not true. Let $v_1, v_2, \cdots, v_{2k+1}$ be the neighbors of $v$ and $v_1, v, v_{2k+1}$ and $u$ be the vertices of the 4-face.
Then the graph $G/ uv_1$ contains a $(2k+1)$-wheel. Hence $G/ uv_1$ is not 3-colorable, a contradiction with Theorem \ref{theorem2.1}. \qed
\end{prof}

\begin{theorem}\label{theorem2.3}
Suppose that $G \in \mathcal{U}_E$ and $G_0$ is a subgraph of $G$. If $G_0$ is uniquely 3-colorable, then we have
\begin{description}
  \item[(i)] $G_0 \in \mathcal{U}_E$;
  \item[(ii)] For any vertex $v\in V(G)\setminus V(G_0)$, $e(v,V(G_0))\leq 2$.
  \end{description}
\end{theorem}
\begin{prof}
 (i) Suppose that $G_0 \notin \mathcal{U}_E$, then there exists an edge $e=uv\in E(G_0)$ such that $G_0-e$ is also uniquely 3-colorable. Let $f$ be a unique 3-coloring of $G$. Since $G \in \mathcal{U}_E$, then $G-e$ has a 3-coloring $f'$ which is distinct from $f$. Note that $f(u)\neq f(v)$, we have $f'(u)\neq f'(v)$.  Thus, $f'$ can be extended to a 3-coloring of $G$. So $G$ has two distinct 3-colorings $f$ and $f'$, which contradicts $G \in \mathcal{U}_E$.

 (ii) Suppose that there exists a vertex $v\in V(G)\setminus V(G_0)$ such that $e(v,V(G_0))=3$. Let $f$ be a unique 3-coloring of $G$ and $v_1,v_2,v_3\in V(G_0)$ be the three neighbors of $v$ in $G$. Then their exist at least two vertices among $v_1,v_2$ and $v_3$ receive the same color. We assume w.l.o.g. that $f(v_1)=f(v_2)$. Since $G \in \mathcal{U}_E$, then $G-vv_1$ has a 3-coloring $f'$ which is distinct from $f$. Note that $G_0$ is uniquely 3-colorable and $f(v)\neq f(v_2)$, we have $f'(v)\neq f'(v_1)$.  Thus, $f'$ can be extended to a 3-coloring of $G$. This is a contradiction. \qed
\end{prof}

\begin{corollary}\label{corollary2.4}
Suppose that $G \in \mathcal{U}_E$ contains a sequence $T_1,T_2,\cdots,T_t$ of triangles satisfying $T_i$ and $T_{i+1}$ have a common edge, where $i=1,2,\cdots, t-1$ and $t\geq2$. Let $v$ and $u$ be the vertices in $V(T_1)\backslash V(T_2)$ and $V(T_t)\backslash V(T_{t-1})$, respectively, then $v\neq u$ and $vu\notin E(G)$.
\end{corollary}
\begin{prof}
Let $v_1,v_2$ be the neighbors of $v$ in $T_1$.
 Since the subgraph of $G$ consists of $t-1$ triangles $T_2,T_3,\cdots,T_{t}$ is uniquely 3-colorable, by Theorem \ref{theorem2.3}, we know that $u$ is not adjacent to $v_1$ or $v_2$ in $G$. Thus, $v\neq u$. Similarly, since the subgraph of $G$ consisting of $t$ triangles $T_1,T_2,\cdots,T_{t}$ is uniquely 3-colorable, we have $vu\notin E(G)$. \qed
\end{prof}

By Corollary \ref{corollary2.4}, we have the following result.

\begin{corollary}\label{corollary2.5}
Suppose that $G \in \mathcal{U}_E$ has no separating 3-cycles. Let $H$ be a subgraph of $G$ that consists of a sequence of triangles $T_1,T_2,\cdots,T_t$ such that each $T_j$ has a common edge with $T_{i}$ for some $i\in\{1,2,\cdots,j-1\}$, where $j=2,3,\cdots, t$. Then $G[V(H)]$ is a maximal outerplanar graph.
\end{corollary}

For a planar graph $G \in \mathcal{U}_E$, if $G$ has no separating 3-cycles, we call the subgraph $H$ in Corollary \ref{corollary2.5} a \emph{triangle-subgraph} of $G$. Note that a triangle is a triangle-subgraph of $G$. Therefore, any $G \in \mathcal{U}_E$ has at least one triangle-subgraph.  A triangle-subgraph $H$ of $G$ is \emph{maximal} if there is no maximal outerplanar subgraph $H'$ of $G$ such that $H\subset H'$. In other words, the graph $H$ consists of the longest sequence $T_1,T_2,\cdots$ of triangles such that each $T_j$ $(j\geq2)$ has a common edge with $T_{i}$ for some $i\in\{1,2,\cdots,j-1\}$.

\begin{theorem}\label{theorem2.6}
Suppose that $G \in \mathcal{U}_E$ has no separating 3-cycles. Let $G_0$ be a uniquely 3-colorable subgraph and $H_1,H_2$ be any two maximal triangle-subgraphs of $G$. If $E(G_0)\cap E(H_i)=\emptyset$, $i=1,2$, then we have
\begin{description}
  \item[(i)] $G_0$ and $H_1$ have at most one common vertex;
  \item[(ii)] If $G_0$ and $H_1$ have a common vertex $v$, then $e(V(G_0-v),V(H_1-v))\leq 1$; otherwise, $e(V(G_0),V(H_1))\leq 3$;
  \item[(iii)] If $H_1$ and $H_2$ have a common vertex $v$, $G_0$ and $H_i$ have a common vertex $v_i$ and $v\neq v_i$, $i=1,2$, then the union of $G_{0},H_{1}$ and $H_{2}$ is uniquely 3-colorable.
  \end{description}
\end{theorem}
\begin{prof} Let $f$ be a unique 3-coloring of $G$.

(i) Suppose, to the contrary, that $G_0$ and $H_1$ have two common vertices $v_1$ and $v_2$. Since $E(G_0)\cap E(H_1)=\emptyset$, then $v_1$ and $v_2$ are not adjacent in both $H_1$ and $G_0$. Otherwise, if $v_1v_2\in E(G_0)\backslash E(H_1)$, this contradicts Corollary \ref{corollary2.4}; if $v_1v_2\in E(H_1)\backslash E(G_0)$, then $G_0+ v_1v_2$ is uniquely 3-colorable but not edge-critical, a contradiction with Theorem \ref{theorem2.3}. By the definition of a triangle-subgraph, we know that there exists a sequence $T_1,T_2,\cdots,T_t$ of triangles in $H_1$ such that $T_i$ and $T_{i+1}$ have a common edge and $\{v_1\}=V(T_1)\backslash V(T_2)$, $\{v_2\}=V(T_t)\backslash V(T_{t-1})$, where $i=1,2,\cdots, t-1$.

 If $f(v_1)=f(v_2)$. Let $v_3$ be a neighbor of $v_1$ in $T_1$, then $v_3\in V(T_2)$. Since $G \in \mathcal{U}_E$, then $G-v_1v_3$ has a 3-coloring $f'$ which is distinct from $f$. Note that both $G_0$ and the subgraph of $H_1$ consists of $t-1$ triangles $T_2,\cdots,T_t$ are uniquely 3-colorable and $f(v_2)\neq f(v_3)$. So $f'(v_1)= f'(v_2)$, $f'(v_2)\neq f'(v_3)$, namely $f'(v_1)\neq f'(v_3)$. Therefore, $f'$ can be extended to a 3-coloring $f'$ of $G$ which is distinct from $f$. This contradicts $G \in \mathcal{U}_E$.

 If $f(v_1)\neq f(v_2)$. Let $v_4$ be a neighbor of $v_1$ in $T_1$ satisfying $f(v_4)=f(v_2)$. Since $G \in \mathcal{U}_E$, then $G-v_1v_4$ has a 3-coloring $f'$ which is distinct from $f$. Since both $G_0$ and the subgraph of $H_1$ consists of $t-1$ triangles $T_2,\cdots,T_t$ are uniquely 3-colorable, we have  $f'(v_1)\neq f'(v_4)$. Therefore, $f'$ can be extended to a 3-coloring $f'$ of $G$ which is distinct from $f$. It is a contradiction.

(ii)\textbf{Case 1}. $G_0$ and $H_1$ have a common vertex $v$.

Suppose that $e(V(G_0-v),V(H_1-v))\geq 2$ and $u_1v_1, u_2v_2$ are two edges with $u_1, u_2 \in V(G_0-v)$ and $v_1, v_2 \in V(H_1-v)$. If there exists a vertex  $u \in \{u_1, v_1, u_2,v_2\}$ such that $f(v)=f(u)$, we assume w.l.o.g. that $u=u_1$, then $f(v)\neq f(v_1)$. Since $G \in \mathcal{U}_E$, $G-u_1v_1$ has a 3-coloring $f'$ which is distinct from $f$. Note that both $G_0$ and $H_1$ are uniquely 3-colorable, we have $f'(v)=f'(u_1)$ and $f'(v)\neq f'(v_1)$. Thus $f'(u_1)\neq f'(v_1)$ and then $f'$ can be extended to a 3-coloring $f'$ of $G$ which is distinct from $f$. If $f(v)\neq f(w)$ for any $w \in \{u_1, v_1, u_2,v_2\}$, then $\{f(u_1),f(v_1)\}=\{f(u_2),f(v_2)\}$. Thus, we have either $f(u_1)=f(u_2)$ and $f(v_1)=f(v_2)$, or $f(u_1)= f(v_2)$ and $f(u_2)= f(v_1)$. Since $G-u_2v_2$ has a 3-coloring $f'$ which is distinct from $f$, and $G_0$ and $H_1$ are uniquely 3-colorable, we have $f'(u_2)\neq f'(v_2)$. Therefore, $f'$ can be extended to a 3-coloring $f'$ of $G$ which is distinct from $f$.

\textbf{Case 2}. $G_0$ and $H_1$ have no common vertex.

 Suppose that $e(V(G_0),V(H_1))\geq 4$ and $u_1v_1, u_2v_2, u_3v_3, u_4v_4$ are 4 edges with $u_i \in V(G_0)$ and $v_i \in V(H_1)$, $i=1,2,3,4$. Then there exist two edges, say $u_1v_1$ and $u_2v_2$, such that $\{f(u_1),f(v_1)\}=\{f(u_2),f(v_2)\}$. By using a similar argument to Case 1, we can obtain a 3-coloring $f'$ of $G-u_1v_1$, which is distinct from $f$ and can be extended to a 3-coloring of $G$. It is a contradiction.

(iii) By definition of $H_1$, there exists a sequence $T_{1},T_{2},\cdots,T_{t}$ of triangles in $H_{1}$ such that $T_{i}$ and $T_{i+1}$ have a common edge and $\{v\}=V(T_{1})\backslash V(T_{2})$, $\{v_1\}=V(T_t)\backslash V(T_{t-1})$, where $i=1,2,\cdots, t-1$.

Suppose that $|\{f(v),f(v_1)$, $f(v_2)\}|=1$. Let $u$ be an arbitrary neighbor of $v$ in $V(T_{1})$. Then $f(v_1)\neq f(u)$. Since $G \in \mathcal{U}_E$, $G-vu$  has a 3-coloring $f'$ which is distinct from $f$. Note that $G_0,H_{2}$ and the subgraph of $H_1$ consists of $t-1$ triangles $T_2,\cdots,T_t$ are uniquely 3-colorable, we have $f'(v)=f'(v_2)=f'(v_1)$ and $f'(v_1)\neq f'(u)$. Therefore, $f'$ can be extended to a 3-coloring $f'$ of $G$ which is distinct from $f$. This contradicts $G \in \mathcal{U}_E$.

Suppose that $|\{f(v),f(v_1),f(v_2)\}|=2$, then their exists $i\in \{1,2\}$ such that $f(v)\neq f(v_i)$. We assume w.l.o.g. that $f(v)\neq f(v_1)$. Let $u$ be a neighbor of $v$ in $V(T_{1})$ satisfying $f(u)=f(v_1)$. If $f(v)=f(v_2)$, then  $f(v_1)\neq f(v_2)$. Since $G \in \mathcal{U}_E$, $G-vu$  has a 3-coloring $f'$ which is distinct from $f$. Note that $G_0,H_{2}$ and the subgraph of $H_1$ consists of $t-1$ triangles $T_2,\cdots,T_t$ are uniquely 3-colorable, we have $f'(v)=f'(v_2)$, $f'(v_1)=f'(u)$ and $f'(v_1)\neq f'(v_2)$. Thus, $f'(v)\neq f'(u)$. Therefore, $f'$ can be extended to a 3-coloring $f'$ of $G$ which is distinct from $f$. This contradicts $G \in \mathcal{U}_E$. If$f(v)\neq f(v_2)$, then $f(v_1)=f(v_2)$. Since $G \in \mathcal{U}_E$, $G-vu$  has a 3-coloring $f'$ which is distinct from $f$. Since $G_0,H_{2}$ and the subgraph of $H_1$ consists of $t-1$ triangles $T_2,\cdots,T_t$ are uniquely 3-colorable, we have $f'(u)=f'(v_1)=f'(v_2)$ and $f'(v)\neq f'(v_2)$. Thus, $f'(v)\neq f'(u)$. This contradicts $G \in \mathcal{U}_E$.

Suppose that $|\{f(v_1),f(v_2),f(v_3)\}|=3$. Using the fact that any coloring $f'$ of two vertices $u,w \in V(G')$ with $f'(u)\neq f'(w)$ can be extended uniquely to a 3-coloring of $G'$, we can obtain that the union of $G_0,H_{1}$ and $H_{2}$ is uniquely 3-colorable, where $G'\in \{G_0,H_1,H_2\}$.
 \qed
\end{prof}

\section{Size of edge-critical uniquely $3$-colorable planar graphs}

In this section, we consider the upper bound of $size(n)$ for edge-critical uniquely 3-colorable planar graphs with $n(\geq 6)$ vertices.

Suppose that $G \in \mathcal{U}_E$ and $G$ has no separating 3-cycles. Let $H_1,H_2,\cdots,H_k$ be all of the maximal triangle-subgraphs of $G$. For two maximal triangle-subgraphs $H_{i}$ and $H_{j}$ having a common vertex $v$, if there exists $H_{\ell}$ such that $H_{i},H_{j}$ and $H_{\ell}$ satisfy the condition of Case (iii) in Theorem \ref{theorem2.6}, namely $H_{i}$ and $H_{\ell}$ have a common vertex (say $v_i$), $H_{j}$ and $H_{\ell}$ have a common vertex (say $v_j$) and $ v_i\neq v \neq v_{j}$, then we say that $H_{i}$ and $H_{j}$ \emph{satisfy} Property \textbf{P}.
Let $G'=H_1\cup H_2\cup \cdots\cup H_k$. (We will use such notation without mention in what follows.)
Now we analyse the relationship between $|F_{\geq 4}(G')|$, the number of $\geq$4-faces of $G'$, and $k$. For a vertex $u\in V(G')$, we use $D(u)$ to denote the number of maximal triangle-subgraphs of $G$ that contain $u$.

\begin{figure}[H]
\begin{center}
  \centering
  \includegraphics[width=320pt]{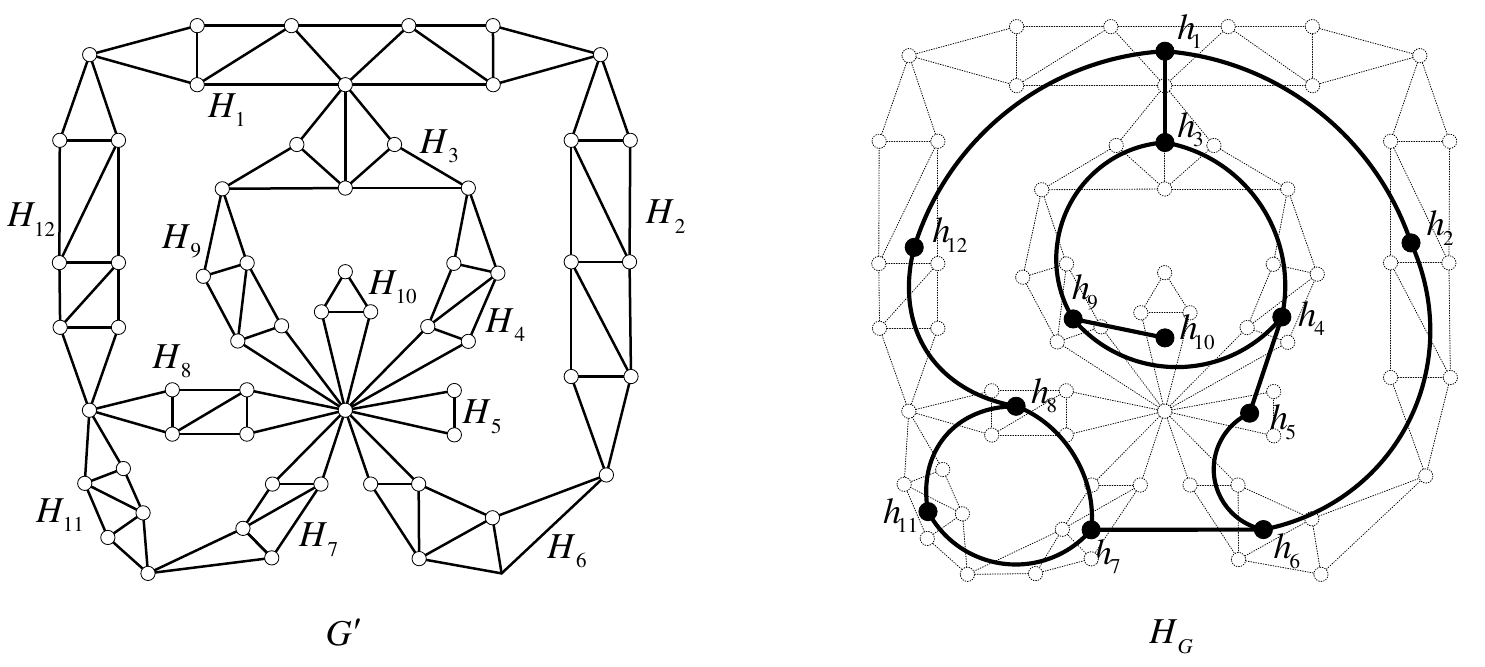}\\
 % \caption{}\label{}
\end{center}
\caption{An example of a graph $G'$ and the corresponding graph $H_G$.}\label{figure1}
\end{figure}

 First we construct a new graph $H_G$ from $G'$ with
$V(H_G)=\{h_1,h_2$, $\cdots$, $h_k\}$, where $h_i$ in $H_G$ corresponds to $H_i$ in $G'$ for any $i \in\{1,2,\cdots,k\}$. The edges in $H_G$ are constructed by the following two steps. \\
\textbf{ Step 1:} For every $u\in V(G')$ with $D(u)=2$, add the edge $h_{i_1}h_{i_2}$ to $H_G$ if both $H_{i_1}$ and $H_{i_2}$ contain $u$. (see e.g. Fig. \ref{figure1}) \\
\textbf{ Step 2:} For every $u\in V(G')$ with $D(u)\geq 3$, let $H_{i_1},H_{i_2}, \cdots, H_{i_{D(u)}}$ contain $u$ and they appear in clockwise order around $u$.
For any $1\leq j< k\leq D(u)$ with $H_{i_j}$ and $H_{i_k}$ satisfying Property \textbf{P}, then add  the edge $h_{i_j}h_{i_k}$ to $H_G$.

Let $G_u$ be the subgraph of $H_G$ with vertex set $V(G_u)=\{h_{i_1},h_{i_2}, \cdots, h_{i_{D(u)}}\}$ and edge set
$E(G_u)=\{h_{i_j}h_{i_k}: \textrm{$H_{i_j}$ and $H_{i_k}$ satisfy Property \textbf{P}}, 1\leq j< k\leq D(u)\}$.

Then we add some edges in $\{h_{i_\ell}h_{i_{\ell+1}}: \ell=1,2,\cdots, D(u)\}$ to $G_u$ such that the resulting graph, denoted by $G_{\langle u \rangle}$, is connected and has the minimum number of edges.
Now the construction of the edges of the graph $H_G$ is completed.
(see e.g. in Fig. \ref{figure1}, we first join the edges $h_4h_9$, $h_7h_8$ and $h_8h_{11}$, then join the edges $h_4h_5$, $h_5h_6$, $h_6h_7$, $h_9h_{10}$ and $h_8h_{12}$.)

\textbf{Remark}. By the definition of $H_G$, if $h_{i}h_{j}\in E(H_G)$, then $H_{i}$ and $H_{j}$ must have a common vertex. For a edge-critical uniquely 3-colorable graph $G$, if $D(u)\leq 2$ for any $u\in G'$, then the graph $H_G$ obtained by above construction is unique; otherwise, $H_G$ is not unique. Furthermore, we have  Theorem \ref{theorem3.1}.

\begin{theorem}\label{theorem3.1}
Suppose that $G \in \mathcal{U}_E$ has no separating 3-cycles. Let $H_1,H_2$, $\cdots$, $H_k$ be all of the maximal triangle-subgraphs of $G$, then $H_G$ is a simple planar graph and $|F(H_G)|=|F_{\geq 4}(G')|$.
\end{theorem}
\begin{prof}
 By Corollary \ref{corollary2.5} and Theorem \ref{theorem2.6}(i), we know that $H_G$ has no loops or parallel edges. So $H_G$ is a simple graph. Note that $G'=H_1\cup H_2\cup \cdots\cup H_k$ is a planar graph.  For any $u\in V(G')$ with $D(u)\geq 3$ and any $1\leq j< k\leq D(u)$ with $H_{i_j}$ and $H_{i_k}$ satisfying Property \textbf{P}, then $G_u$ is a planar graph and there exist no edges $h_{i_a}h_{i_b},h_{i_c}h_{i_d}\in E(G_u)$ such that $i_a \in \{i_c+1,\cdots,i_d-1\}$ and $i_b \in \{i_d+1,\cdots,i_c-1\}$, where $a,b,c,d\in \{1,2,\cdots,D(u)\}$ and the subscripts are taken modulo $D(u)$. Now we prove that $G_u$ is a forest. If $h_{i'_1}h_{i'_2},h_{i'_1}h_{i'_3}\in E(G_u)$, then, by the definition of $G_u$ and Theorem \ref{theorem2.6}(iii), we know that there exist $H_{\ell_1}$ and $H_{\ell_2}$ such that the graph $H_{i'_1}\cup H_{i'_2}\cup H_{i'_3}\cup H_{\ell_1}\cup H_{\ell_2}$ is uniquely 3-colorable. Thus, $H_{i'_2}$ and $H_{i'_3}$ does not satisfy Property \textbf{P}, namely $h_{i'_2}h_{i'_3}\notin E(G_u)$. Therefore, $G_u$ is a forest. By the definition of $G_{\langle u \rangle}$, it is easy to see that $G_{\langle u \rangle}$ is a tree. By the definition of $H_G$, we can conclude that $H_G$ is a planar graph.

For any distinct faces $f_1$ and $f_2$ of $H_G$, by the definition of $H_G$, it can be seen that there exist two distinct $\geq$4-faces of $G'$ corresponding to $f_1$ and $f_2$, respectively.
Conversely, for any $\geq$4-face $f'$ of $G'$, let $H_{j_1},H_{j_2}, \cdots, H_{j_{t}}$ be all of the maximal triangle-subgraphs satisfying $H_{j_\ell}$ and $f'$ have common edges, $\ell=1,2,\cdots, t$. Let $u_\ell$ be the common vertex of $H_{j_\ell},H_{j_{\ell+1}}$, because $G_{\langle u_\ell\rangle}$ is tree, there exists a unique face of $H_G$ incident with $h_{j_1},h_{j_2}, \cdots, h_{j_{t}}$. Thus, $|F(H_G)|=|F_{\geq 4}(G')|$.
\qed
\end{prof}

\begin{theorem}\label{theorem3.2}
Suppose that $G \in \mathcal{U}_E$ has no separating 3-cycles. Let $f_0,f_1$, $\cdots$, $f_t$ be a sequence of faces in $H_G$ such that $f_\ell$ and $f_{m}$ are adjacent, $\ell=1,2,\cdots, t$, $m\in \{0,1,\cdots,\ell-1\}$. If $d(f_0)=3$ and $d(f_\ell)=4$, $\ell=1,2,\cdots, t$, let $h_{i_1},h_{i_2},\cdots, h_{i_s}$ be all of the vertices incident with the faces $f_0,f_1,\cdots,f_{t}$. Then $|V(f_t)\setminus \bigcup_{\ell=0}^{t-1}V(f_{\ell})|=2$ and $H_{i_1}\cup H_{i_2}\cup\cdots\cup H_{i_s}$ is uniquely 3-colorable, where $V(f)$ denotes the set of the vertices incident with $f$.
\end{theorem}
\begin{prof}
The proof is by induction on $t$. Let $V(f_0)=\{h_{i_1},h_{i_2},h_{i_3}\}$ and $h_{j_{s-r}},\cdots, h_{i_s}$ be the vertices incident with $f_{t}$, but not incident with $f_0,f_1$, $\cdots$, $f_{t-1}$. If $t=0$, since $d(f_0)=3$, by Theorem \ref{theorem2.6} (iii), we know that $H_{i_1}\cup H_{i_2}\cup H_{i_3}$ is uniquely 3-colorable.
Suppose that $t\geq 1$. By hypothesis, $H_{i_1}\cup H_{i_2}\cup\cdots\cup H_{i_{s-r-1}}$ is uniquely 3-colorable.  Since $d(f_t)=4$, by Theorem \ref{theorem2.6} (i), we have $r=1$, namely $|V(f_t)\setminus \bigcup_{\ell=0}^{t-1} V(f_{\ell})|=2$. Therefore, by Theorem \ref{theorem2.6} (iii), we obtain that $H_{i_1}\cup H_{i_2}\cup\cdots\cup H_{i_s}$ is uniquely 3-colorable.
\qed
\end{prof}

\vspace{0.3cm}

For a planar graph $G$, let $C$ and $C'$ be two cycles of $G$.  $C$ and $C'$ are \emph{dependent} if there exists a sequence $C_1(=C),C_2,\cdots,C_t(=C')$ of cycles of $G$ such that $C_\ell$ and $C_{\ell+1}$ have common edges and $|V(C_s)|=4$, where $\ell=1,2,\cdots, t-1$, $s=2,3,\cdots, t-1$. Obviously, if $C$ and $C'$ have a common edge, then they are dependent.

\begin{lemma}\label{lemma3.3}
Let $G$ be a planar graph, $|V(G)|\geq 4$. If any $i$-cycle of $G$ is dependent with at most $i-3$ 3-cycles for $3\leq i \leq5$ and with at most $i-2$ 3-cycles for $i\geq 6$, then $|V(G)|\geq |F(G)|+2$.
\end{lemma}
\begin{prof}
The proof is by contradiction. Let $G$ be a smallest counterexample to the lemma, then $G$ satisfies the conditions of the lemma and $|V(G)|< |F(G)|+2$. Suppose that $G$ is not connected, let $G_1$ be a connected component of $G$. If $|V(G_1)|\leq 3$ and $|V(G-V(G_1))|\leq 3$, it is easy to see that $|V(G)|\geq |F(G)|+2$. This is a contradiction. Otherwise, we assume w.l.o.g. that $|V(G-V(G_1))|\geq 4$. Since any $i$-cycle of $G$ is dependent with at most $i-3$ 3-cycles for $3\leq i \leq5$ and with at most $i-2$ 3-cycles for $i\geq 6$, the same is true of $G_1$ and $G-V(G_1)$. By the minimality of $G$, we have $|V(G-V(G_1))|\geq |F(G-V(G_1))|+2$. Furthermore, if $|V(G_1)|\geq 4$, then $|V(G_1)|\geq |F(G_1)|+2$; otherwise, $|V(G_1)|\geq |F(G_1)|$. Therefore, $|V(G)|=|V(G_1)|+|V(G-V(G_1))|\geq |F(G_1)|+|F(G-V(G_1))|+2= |F(G)|+3$, a contradiction.

Suppose that $G$ is connected. If $G$ contains a cut vertex $u$, let $V_1,V_2,\cdots,V_r$ be the vertex sets of the connected components of $G-u$, respectively, and $G_j=G[\{u\}\cup V_i]$, $j=1,2,\cdots,r$. Obviously, $G_j$ satisfies the conditions of the lemma. If $|V(G_j)|\geq 4$, then, by the minimality of $G$, $|V(G_j)|\geq |F(G_j)|+2$; otherwise, $|V(G_j)|\geq |F(G_j)|+1$, $j=1,2,\cdots,r$.  Therefore, $|V(G)|=\sum_{j=1}^{r} |V(G_j)|-(r-1)\geq \sum_{j=1}^{r} (|F(G_j)|+1)-(r-1)= |F(G)|+(r-1)+1\geq |F(G)|+2$, a contradiction.

Now we assume that $G$ is 2-connected. If $G$ contains no 3-faces, then $2|E(G)|=\sum_{f\in F(G)} d(f)\geq 4|F(G)|$. Thus $|E(G)|\geq 2|F(G)|$. By Euler's Formula, we have $|V(G)|\geq |F(G)|+2$. This contradicts the choice of $G$. If $G$ contains exactly one 3-face, then $G$ contains at least one $\geq$5-face. Thus, $2|E(G)|=\sum_{f\in F(G)} d(f)\geq 4|F(G)|$ and then $|V(G)|\geq |F(G)|+2$.  If $G$ contains at least two 3-faces, then each 3-face is dependent with at least two $\geq$5-faces for $G$ is 2-connected. We claim that $|E(G)|\geq 2|F(G)|$, namely $\sum_{f\in F(G)}d(f)-4\geq 0$.

 For any face $f\in F(G)$, we set the \emph{initial charge} of $f$ to be $ch(f)=d(f)-4$. We now use the discharging procedure, leading to the final charge $ch'$, defined by applying the following rule:

 \textbf{RULE}. Each 3-face receives $\frac{1}{2}$ from each dependent $\geq$5-face.

For any face $f\in F(G)$, if $d(f)=3$, since $f$ is dependent with at least two $\geq$5-faces, then $ch'(f)\geq ch(f)+2\times \frac{1}{2}=0$. If $d(f)=4$, then $ch'(f)=ch(f)=0$. If $d(f)=5$, then $ch'(f)\geq ch(f)-\frac{1}{2}\cdot [d(f)-3]=0$. If $d(f)\geq 6$, then, by hypothesis, $ch'(f)\geq ch(f)-\frac{1}{2}\cdot [d(f)-2]\geq \frac{1}{2}\cdot d(f)-3\geq 0$. Therefore, $\sum_{f\in F(G)}ch(f)=\sum_{f\in F(G)}ch'(f)\geq 0$.

Thus, by Euler's Formula, we have $|V(G)|\geq |F(G)|+2$. This contradicts the choice of $G$.
  \qed
\end{prof}

\begin{theorem}\label{theorem3.4}
Suppose that $G \in \mathcal{U}_E$ has no separating 3-cycles. If $G$ has $k$ maximal triangle-subgraphs $H_1,H_2$, $\cdots$, $H_k$ and $k\geq 4$, then $|F(H_G)|\leq |V(H_G)|-2$.
\end{theorem}
\begin{prof}
By Lemma \ref{lemma3.3}, it suffices to prove that any $i$-cycle of $H_G$ is dependent with at most $i-3$ 3-cycles if $3\leq i \leq5$ and with at most $i-2$ 3-cycles if $i\geq 6$. The proof is by contradiction. Let $C$ be a $i$-cycle of $H_G$, and $C$ is dependent with at least $i-2$ 3-cycles ($3\leq i \leq5$) or with at least $i-1$ 3-cycles ($i\geq 6$). If $i=3$ or 4, by Theorems \ref{theorem2.6} (i) and \ref{theorem3.2}, it is easy to see that there exist no dependent 3-cycles and at most one 3-cycle that is dependent with a 4-cycle. This contradicts the hypothesis. Suppose that $i\geq 5$, let $r=i-2$ if $i=5$ and $r=i-1$ if $i\geq 6$. Let $C=C_{j,0},C_{j,1},\cdots,C_{j,t_j}$ be a sequence of cycles of $H_G$ such that $C_{j,\ell}$ and $C_{j,\ell+1}$ have common edges, $|V(C_{j,s})|=4$ and $|V(C_{j,t_j})|=3$, where $\ell=0,1,\cdots, t_j-1$, $s=1,2,\cdots, t_j-1$ and $j=1,2,\cdots, r$. Then for any $a\in \{1,2,\cdots, t_1\}$ and $b\in \{1,2,\cdots, t_2\}$, $C_{1,a}$ and $C_{2,b}$ are not dependent. Otherwise, $C_{1,a}$ is dependent with two 3-cycles $C_{1,t_1}$ and $C_{2,t_2}$. Therefore, each pair of 4-cycles in $\{C_{1,1},C_{2,1},\cdots,C_{r,1}\}$ have no common edges. Moreover, $C$ and $C_{j,1}$ have exactly one common edge because $r\geq i-2$, $j=1,2,\cdots, r$.

 If $i=5$, then $r=3$. We assume w.l.o.g. that $V(C)=\{h_1,h_2,\cdots,h_5\}$ and let $V_j=\bigcup_{\ell=1}^{t_j}V(C_{j,\ell})$. Now we consider the following two cases:

\textbf{Case 1}. $C_{1,1}\cup C_{2,1}\cup C_{3,1}$ contains 4 vertices of $C$, see Fig.\ref{figure2} (a). Assume w.l.o.g. that $h_1\notin V(C_{1,1}\cup C_{2,1}\cup C_{3,1})$ and $V_1\cup V_2 \cup V_3=\{h_2,h_3,\cdots,h_p\}$. By Theorem \ref{theorem3.2}, we can conclude that $G_0=H_2\cup H_3\cup\cdots\cup H_p$ is uniquely 3-colorable. Note that $H_1$ and $G_0$ have two common vertices, this contradicts Theorem \ref{theorem2.6} (i).

\textbf{Case 2}. $C_{1,1}\cup C_{2,1}\cup C_{3,1}$ contains 5 vertices of $C$, see Fig.\ref{figure2} (b).  Assume w.l.o.g. that $V_2\cup V_3=\{h_3,h_4,\cdots,h_{p}\}$, $V_1=\{h_1,h_2,h_{p+1},\cdots,h_{p'}\}$ and $h_{p'}\in V(C_{1,t_1})\setminus V(C_{1,t_{1}-1})$.  By Theorem \ref{theorem3.2}, we can conclude that $G_0=H_3\cup H_4\cup\cdots\cup H_p$ is uniquely 3-colorable. Then, by Theorem \ref{theorem2.6} (iii), we obtain that $G_1=G_0\cup H_1\cup H_2\cup H_{p+1}\cup \cdots\cup H_{p'-1}$ is uniquely 3-colorable. Note that $H_p'$ and $G_1$ have two common vertices, this contradicts Theorem \ref{theorem2.6} (i).

\begin{figure}[H]
\begin{center}
  \centering
  \includegraphics[width=320pt]{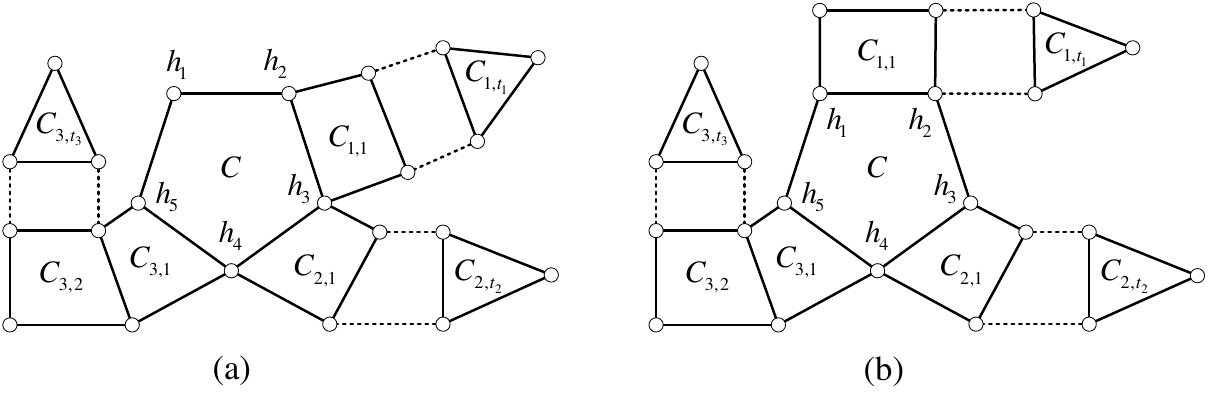}\\
 % \caption{}\label{}
\end{center}
\caption{Two cases of a $5$-cycle dependent with three 3-cycles.}\label{figure2}
\end{figure}

If $i\geq 6$, then $r=i-1$. Assume w.l.o.g. that $V(C)=\{h_1,h_2,\cdots,h_i\}$. In this case, there exists only one edge, say $h_1h_i$, of $C$ that is not in $C_{1,1}\cup C_{2,1}\cup\cdots \cup C_{r,1}$. Suppose that $C_{j,1}$ contains the edge $h_{j}h_{j+1}$ and $V_2\cup V_3 \cup\cdots \cup V_{r-1}=\{h_2,h_3,\cdots,h_p\}$, $j=1,2,\cdots,r-1$. By Theorem \ref{theorem3.2}, we can conclude that $G_0=H_2\cup H_3\cup\cdots\cup H_p$ is uniquely 3-colorable. Note that $H_1$ and $G_0$ have two common vertices, this contradicts Theorem \ref{theorem2.6} (i).   \qed
\end{prof}

By Theorems \ref{theorem3.1} and \ref{theorem3.4}, we obtain the following Corollary \ref{corollary3.5}.

\begin{corollary}\label{corollary3.5}
Suppose that $G \in \mathcal{U}_E$ has no separating 3-cycles. If $G$ has $k$ maximal triangle-subgraphs $H_1,H_2$, $\cdots$, $H_k$ and $k\geq 4$, then $|F_{\geq 4}(G')|\leq k-2$.
\end{corollary}

\begin{theorem}\label{theorem3.6}
Let $G \in \mathcal{U}_E$ and $|V(G)|\geq 6$, then $|E(G)|\leq \frac{5}{2}|V(G)|-6$.
\end{theorem}
\begin{prof}
The proof is by induction on $n=|V(G)|$. It is easy to check that the theorem is true for $n=6$. Suppose that the theorem is true for all
edge-critical uniquely 3-colorable planar graphs with $p$ vertices, where $6\leq p \leq n-1$ and $n\geq 7$. Let $G\in \mathcal{U}_E$ and $|V(G)|=n$.
We consider the following two cases:

\textbf{Case 1}. $G$ contains a separating 3-cycle $C$.

Let $G_1$ (resp. $G_2$) be the subgraph of $G$ consists of $C$ together with its interior (resp. exterior). Then both $G_1$ and $G_2$ are uniquely 3-colorable planar graphs. Otherwise, suppose that $G_1$ has two distinct 3-colorings, then each 3-coloring of $G_1$ can be extended to a 3-coloring of $G$. This contradicts $G \in \mathcal{U}_E$. By Theorem \ref{theorem2.3}, we have $G_1,G_2 \in \mathcal{U}_E$. If $V(G_i)\geq 6$ for $i=1,2$, by induction, $|E(G_i)|\leq \frac{5}{2}|V(G_i)|-6$. Thus, $|E(G)|=|E(G_1)|+|E(G_2)|-3\leq \frac{5}{2}|V(G_1)|-6+\frac{5}{2}|V(G_2)|-6-3\leq \frac{5}{2}(|V(G)|+3)-15< \frac{5}{2}|V(G)|-6$. If $V(G_1)\leq 5$ or $V(G_2)\leq 5$, since $G_1,G_2 \in \mathcal{U}_E$, there exists a vertex $v$ in $V(G_1)\setminus V(C)$ or $V(G_2)\setminus V(C)$ such that $d_G(v)=2$. Therefore, $G- v$ is uniquely 3-colorable and then $G- v \in \mathcal{U}_E$. By induction, $|E(G- v)|\leq \frac{5}{2}|V(G- v)|-6$. Thus, $|E(G)|=|E(G- v)|+2\leq \frac{5}{2}(|V(G)|-1)-4< \frac{5}{2}|V(G)|-6$.

\textbf{Case 2}. $G$ contains no separating 3-cycles.

Using the fact that every planar graph with $n$ vertices is a subgraph of a maximal planar graph with the same vertices, we may assume that $G_{max}$ is a maximal planar graph with $n$ vertices and $G$ is a subgraph of $G_{max}$. Let $q=|E(G_{max})|-|E(G)|$, then $|E(G)|=3n-6-q$ and $|F(G)|=2n-4-q$. In this case, we prove the theorem by showing that $q\geq \frac{n}{2}$.

Let $H_1,H_2$, $\cdots$, $H_k$ be all of the maximal triangle-subgraphs of $G$, $G'=H_1\cup H_2\cup \cdots\cup H_k$ and $H_i$ contain $t_i$ 3-faces, where $i=1,2\cdots,k$. Then $|V(H_i)|=t_i+2$, $|E(H_i)|=2t_i+1$ and $|F_3(G)|=\sum_{i=1}^{k}t_i$. Moreover, $|E(G')|=\sum_{i=1}^{k}|E(H_i)|=\sum_{i=1}^{k}(2t_i+1)=k+2|F_3(G)|$.

Let $G^*$ be the dual graph of $G$ and $G^*_0$ be the subgraph of $G^*$ induced by $V_{\geq4}(G^*)$, the set of vertices of degree at least 4 in $G^*$. By Euler's Formula, we have $|E(G^*_0)|=|V(G^*_0)|+|F(G^*_0)|-\omega(G^*_0)-1$, where $\omega(G^*_0)$ is the number of connected components of $G^*_0$. By the definition of $G^*$, we have $|V(G^*_0)|=|F_{\geq4}(G)|$, $|E(G^*_0)|=|E(G)|-|E(G')|$ and $|F(G^*_0)|=|V(G)|-|V(G')|+\omega(G')$. Since $G$ contains no separating 3-cycles, then $\omega(G^*_0)=|F_{\geq 4}(G')|$. Therefore,

\begin{equation}\label{equ1}
\begin{split}
         & |E(G)| = |E(G')|+|E(G_0^*)|\\
         =& k+ 2|F_3(G)|+|F_{\geq4}(G)|+|V(G)|-|V(G')|+\omega(G')-|F_{\geq 4}(G')|-1\\
         =&2n-4-q +|F_3(G)|+(k-|F_{\geq 4}(G')|)+(n-|V(G')|)+(\omega(G')-1)
\end{split}
\end{equation}

Note that $n-|V(G')|\geq 0$ and $\omega(G')-1\geq 0$ in Formula (\ref{equ1}).
Because $G_{max}$ has $2n-4$ 3-faces by Euler's Formula and removing a edge decreases the number of 3-faces by at most two, we have

\begin{equation}\label{equ2}
        |F_3(G)|\geq 2n-4-2q
\end{equation}

Suppose that $k=1$, then $|F_{\geq 4}(G')|=\omega(G')=1$, $H_1$ is a maximal outerplanar graph and $H_1=G[V(H_1)]$ by Corollary \ref{corollary2.5}. If $|V(G')|=n$, then $G=H_1$. In this case, $|E(G)|=2n-3<\frac{5}{2}n-6$ since $n\geq 7$. If $|V(G')|=n-1$, then, by Theorem \ref{theorem2.3}, $|E(G)|=|E(H_1)|+2=2(n-1)-3+2<\frac{5}{2}n-6$. If $|V(G')|\leq n-2$, then, by Formula (\ref{equ1}), we have $|E(G)|\geq 2n-4-q +2n-4-2q+2=4n-3q-6$. Since $|E(G)|=3n-6-q$, we have $q\geq \frac{n}{2}$. Therefore, $|E(G)|\leq \frac{5}{2}n-6$.

Suppose that $k=2$, by Theorem \ref{theorem2.6}(i), we have $|F_{\geq 4}(G')|=1$ and $\omega(G')\leq 2$. If $\omega(G')=1$ and $|V(G')|=n$, then $H_1$ and $H_2$ have a common vertex. By  Theorem \ref{theorem2.6}(ii), there exists at most one edge in $E(G)\setminus (E(H_1)\cup E(H_2))$. Therefore, $|E(G)|\leq |E(H_1)|+|E(H_2)|+1=2t_1-3+2t_2-3 +1=2(n+1)-5<\frac{5}{2}n-6$. If $\omega(G')=2$ or $|V(G')|\leq n-1$,  then, by Formula (\ref{equ1}), we have $|E(G)|\geq 2n-4-q +2n-4-2q+1+1=4n-3q-6$. Similarly, we can obtain $q\geq \frac{n}{2}$, and hence, $|E(G)|\leq \frac{5}{2}n-6$.

Suppose that $k=3$, by Theorem \ref{theorem2.6}(i) and (iii), we have $|F_{\geq 4}(G')|\leq 2$ and $\omega(G')\leq 3$. If  $|F_{\geq 4}(G')|=1$, then, by Formula (\ref{equ1}), we have $|E(G)|\geq 2n-4-q +2n-4-2q+2=4n-3q-6$. Therefore, $q\geq \frac{n}{2}$ and then $|E(G)|\leq \frac{5}{2}n-6$. If  $|F_{\geq 4}(G')|=2$, then, by Theorem \ref{theorem2.6}(iii), we know that $G'$ is uniquely 3-colorable. In this case, if $|V(G')|=n$, then $G=G'$ and $|E(G)|=|E(H_1)|+|E(H_2)|+|E(H_3)|=2t_1-3+2t_2-3+2t_3-3=2(n+3)-9<\frac{5}{2}n-6$. If $|V(G')|\leq n-1$, then, by Formula (\ref{equ1}), we have $|E(G)|\geq 2n-4-q +2n-4-2q+1+1=4n-3q-6$. Therefore, $q\geq \frac{n}{2}$ and then $|E(G)|\leq \frac{5}{2}n-6$.

Suppose that $k\geq 4$, by Corollary \ref{corollary3.5}, we have $k-|F_{\geq 4}(G')|\geq 2$. By Formula (\ref{equ1}), we have $|E(G)|\geq 2n-4-q +2n-4-2q+2=4n-3q-6$. Therefore, $q\geq \frac{n}{2}$ and  then $|E(G)|\leq \frac{5}{2}n-6$.
\qed
\end{prof}

\section{Concluding Remarks}

In this section we give some edge-critical uniquely $3$-colorable planar graphs which have $n(=10,12, 14)$ vertices and $\frac{5}{2}n-7$ edges.

\begin{figure}[H]
\begin{center}
  \centering
  \includegraphics[width=80pt]{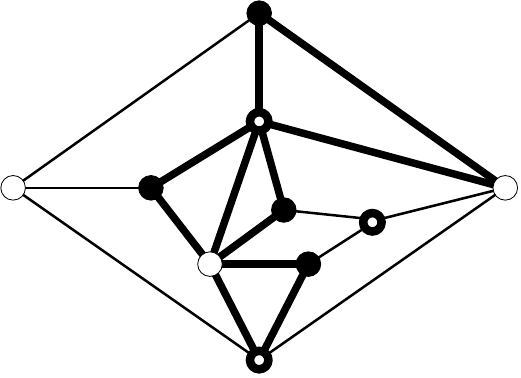}\\
 % \caption{}\label{}
\end{center}
\caption{A edge-critical uniquely $3$-colorable planar graph $G_1$.}\label{figure3}
\end{figure}

Fig. \ref{figure3} shows a edge-critical uniquely $3$-colorable planar graph $G_1$, which has $10$ vertices and 18 edges, and a unique 3-coloring of $G_1$.

\begin{figure}[H]
\begin{center}
  \centering
  \includegraphics[width=200pt]{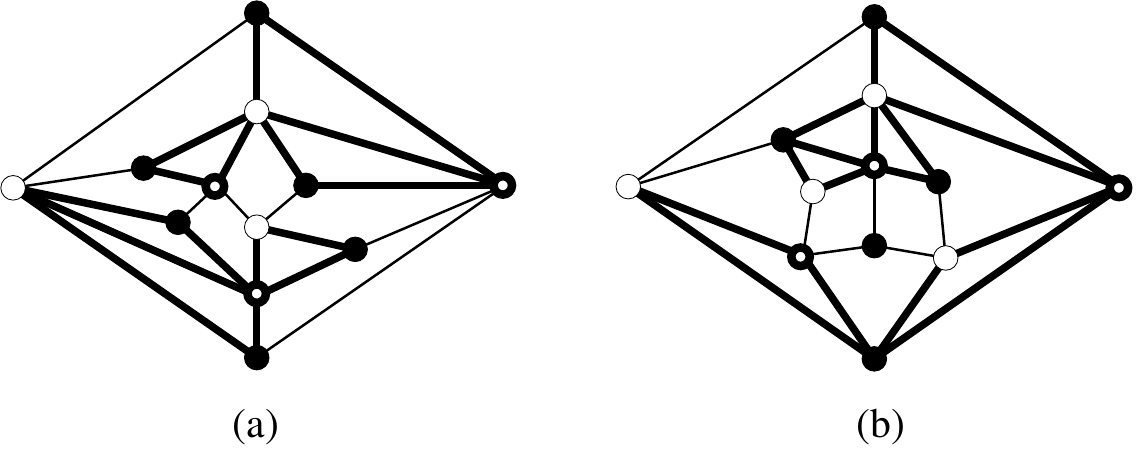}
 % \caption{}\label{}
\end{center}
\caption{Two edge-critical uniquely $3$-colorable planar graphs $G_2$ and $G_3$.}\label{figure4}
\end{figure}

Fig. \ref{figure4} shows two edge-critical uniquely $3$-colorable planar graphs $G_2$ and $G_3$, both of which have 12 vertices and 23 edges, and their unique 3-colorings.

\begin{figure}[H]
\begin{center}
  \centering
  \includegraphics[width=220pt]{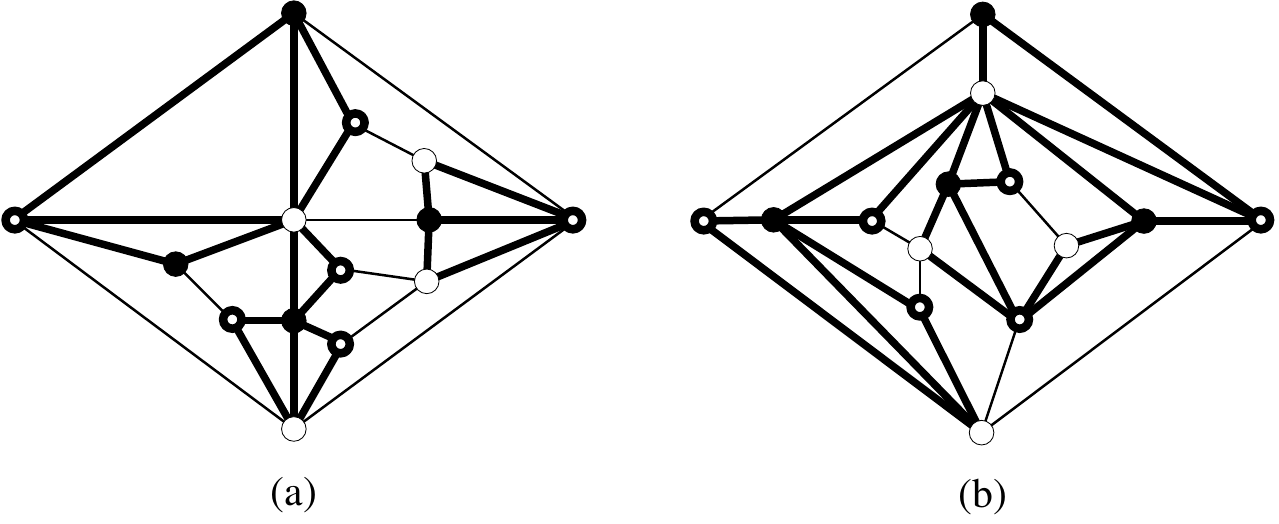}
 % \caption{}\label{}
\end{center}
\caption{Two edge-critical uniquely $3$-colorable planar graphs $G_4$ and $G_5$.}\label{figure5}
\end{figure}

Fig. \ref{figure5} shows two edge-critical uniquely $3$-colorable planar graphs $G_4$ and $G_5$, both of which have $14$ vertices and 28 edges, and their unique 3-colorings.

Note that for $G_1$, we have $k(G_1)-|F_{\geq 4}(G'_1)|=2$, where $k(G_1)$ is the number of maximal triangle-subgraphs of $G_1$. For $i\in \{2,4,5\}$, $|V(G'_i)|=|V(G_i)|$; For $i\in \{2,4\}$, $\omega(G'_i)=1$. Furthermore, we have $|F_3(G_i)|\geq 2|V(G_i)|-4-2q$ for $i\in \{1,2,3,4,5\}$, namely the equality of Formula \ref{equ2} holds for $G_i$.

%It is not hard to check that $size(n)=2n-3$ by Theorem \ref{theorem2.6} if $n\leq 9$. Therefore, by Theorem \ref{theorem3.6}, we have the following corollary.

%\begin{corollary}\label{corollary4.1}
%If $n\geq 10$, then we have
% $$ \frac{5}{2}n-7 \leq size(n)\leq \frac{5}{2}n-6.$$
%\end{corollary}

%\section{Final remarks}

\end{document}